\documentclass[12pt]{amsart}
\usepackage{amsmath,amsthm,amssymb}
\usepackage{graphicx}
\begin{document}
\newtheorem{thm}{Theorem}
\newtheorem{pro}[thm]{Proposition}
\newtheorem{cor}[thm]{Corollary}
\newtheorem{lem}[thm]{Lemma}
\newtheorem{dfn}[thm]{Definition}
\newtheorem{rem}[thm]{Remark}
\newtheorem{prob}[thm]{Problem}
\newtheorem{exam}[thm]{Example}
\newtheorem{conj}[thm]{Conjecture}
\renewcommand{\theequation}{\arabic{section}.\arabic{equation}}
\renewcommand{\labelenumi}{\rm{(\arabic{enumi})}}
\title[Lagrangian Floer homology of a pair of real forms]
{Lagrangian Floer homology of a pair of real forms in
Hermitian symmetric spaces of compact type}
\author[H. Iriyeh, T. Sakai and H. Tasaki]
{Hiroshi Iriyeh, Takashi Sakai and Hiroyuki Tasaki}
\date{}

\begin{abstract}
In this paper we calculate the Lagrangian Floer homology
$HF(L_0, L_1 : {\mathbb Z}_2)$ of a pair of real forms $(L_0,L_1)$
in a monotone Hermitian symmetric space $M$ of compact type
in the case where $L_0$ is not necessarily congruent to $L_1$.
In particular, we have a generalization of the Arnold-Givental inequality
in the case where $M$ is irreducible.
As its application, we prove that the totally geodesic Lagrangian
sphere in the complex hyperquadric is globally volume minimizing under
Hamiltonian deformations.
\end{abstract}

\maketitle

\section{Introduction and main results}

Let $(M,\omega)$ be a symplectic manifold, i.e.,
$M$ is a smooth manifold with a closed nondegenerate $2$-form $\omega$.
Let $L$ be a Lagrangian submanifold in $M$, i.e.,
$\dim_{\mathbb R}L=\frac{1}{2}\dim_{\mathbb R}M$ and
$\omega$ vanishes on $L$.
For a pair of closed Lagrangian submanifolds $(L_0,L_1)$ in $M$,
we can define Lagrangian Floer homology $HF(L_0, L_1 : {\mathbb Z}_2)$
with coefficient ${\mathbb Z}_2$ under some appropriate topological
conditions.

In 1988, Floer \cite{Floer88} defined the homology when
$\pi_2(M,L_i)=0,\ i=0,1$, and proved that it is isomorphic to
the singular homology group $H_*(L_0,{\mathbb Z}_2)$ of $L_0$
in the case where $L_0$ is {\it Hamiltonian isotopic} to $L_1$.
As a result, he solved affirmatively the so called
Arnold conjecture for Lagrangian intersections in that case
(see \cite{Arnold65} and \cite{Floer88}).
A symplectic diffeomorphism $\phi$ of $(M,\omega)$ is called
{\it Hamiltonian} if $\phi$ is represented by the time-$1$ map of
the flow $\{ \phi_t \}$ of a time dependent Hamiltonian vector field on $M$,
i.e., $\frac{d}{dt}\phi_t(x)=X_{H_t}(\phi_t(x)),\ \phi_0(x)=x$,
where $X_{H_t}$ is defined by the equation
$\omega(X_{H_t},\cdot)=dH_t$
for a smooth function
$H:[0,1] \times M \to \mathbb R$.
We denote by $\mathrm{Ham}(M,\omega)$ the set of all Hamiltonian
diffeomorphisms of $M$.

After that, Givental \cite{Givental88} and Chang-Jiang \cite{Chang-Jiang90}
proved the conjecture for $L={\mathbb R}P^n \subset M={\mathbb C}P^n$
independently. (See also \cite{Oh93'}).
In the same paper, Givental posed the following conjecture
which generalizes the above results by Floer and himself.

\begin{conj}[\bf Arnold-Givental] \label{conj:AG} \rm
Let $(M,\omega)$ be a symplectic manifold and
$\tau:M \to M$ be an anti-symplectic involution of $M$.
Assume that the fixed point set $L=\mathrm{Fix}(\tau)$ is not empty
and compact.
Then for any $\phi \in \mathrm{Ham}(M,\omega)$
such that the Lagrangian submanifold $L$ and its image $\phi L$
intersect transversally, the inequality
\begin{equation*}
\#(L \cap \phi L) \geq SB(L, \mathbb Z_2)
\end{equation*}
holds, where $SB(L, \mathbb Z_2)$ denotes the sum of
$\mathbb Z_2$-Betti numbers of $L$.
\hfill\qed
\end{conj}

Note that the assumption of Conjecture \ref{conj:AG}
admits many explicit examples.
For instance,
any real form $L$ of Hermitian symmetric spaces of compact type is included.

The first substantial progress towards Conjecture \ref{conj:AG}
was made by Y.-G. Oh \cite{Oh95}.
He solved the Arnold-Givental conjecture affirmatively for
{\it monotone} real forms of Hermitian symmetric spaces of compact type
(see Corollary \ref{cor:Oh} below).
To prove it, he improved Floer's construction
so as to apply the Lagrangian Floer homology theory to the case of
monotone Lagrangian submanifolds
(see \cite{Oh93} and Section 2 in this paper).
After that, Frauenfelder \cite{Frauenfelder04} proved
the Arnold-Givental conjecture for some class of Lagrangian submanifolds
in Marsden-Weinstein quotients, which are fixed point sets of
some anti-symplectic involution.
Recently, Fukaya, Oh, Ohta and Ono \cite{FOOO} proved
Conjecture \ref{conj:AG} in a considerably more general setting,
but the general case is still an open problem.

For a pair of Lagrangian submanifolds $(L_0,L_1)$,
where $L_0$ is not Hamiltonian isotopic to $L_1$
in a symplectic manifold $(M,\omega)$,
there are relatively few examples where it is known how to calculate
$HF(L_0,L_1:{\mathbb Z}_2)$.
Recently, explicit calculations of the Floer homology of Lagrangian
submanifolds in toric Fano manifolds have been intensively studied
(see \cite{Al08}, \cite{Alston-Amorim10} and \cite{FOOO10}).

In this paper, we shall focus on real forms $L_0,L_1$
of a Hermitian symmetric space of compact type and calculate
the Lagrangian Floer homology $HF(L_0, L_1 : {\mathbb Z}_2)$
in a unified method.
Let $(M,J_0,\omega)$ be a Hermitian symmetric space of compact type.
A submanifold $L$ of $M$ is called a {\it real form} if there exists
an anti-holomorphic involutive isometry $\sigma$ of $M$ satisfying
$$
L = \{x \in M \mid \sigma(x) = x\}.
$$
Note that a real form of a Hermitian symmetric space is
a totally geodesic Lagrangian submanifold.
We denote by $I_0(M)$ the identity component of
the holomorphic isometry group of $M$.
For two subsets $A,B$ of $M$, we say that $A$ is {\it congruent} to $B$
if there exists $g \in I_0(M)$ such that $B=gA$.
Let $L=\mathrm{Fix}(\sigma)$ be a real form and
$g$ a holomorphic isometry of $M$.
Then $gL=\mathrm{Fix}(g \sigma g^{-1})$ is also a real form of $M$.

The following is the main result.

\begin{thm} \label{thm:main1}
Let $(M,J_0,\omega)$ be a Hermitian symmetric space of compact type
which is monotone as a symplectic manifold.
Let $L_0, L_1$ be real forms of $M$ such that $L_0$ intersects $L_1$
transversally.
Assume that the minimal Maslov numbers of $L_0$ and $L_1$ are greater
than or equal to $3$. Then we have
\begin{equation*}
HF(L_0, L_1 : {\mathbb Z}_2)
\cong \bigoplus_{p \in L_0 \cap L_1} {\mathbb Z}_2 [p].
\end{equation*}
That is, the intersection $L_0 \cap L_1$ itself becomes
a basis of the Floer homology $HF(L_0, L_1 : {\mathbb Z}_2)$.
\end{thm}

If $M$ is irreducible, then the assumptions for $M$, $L_0$ and $L_1$
are satisfied automatically except for only one case (see Section 3).
Moreover, using the structure of the transverse intersection
$L_0 \cap L_1$ which was examined by Tanaka and Tasaki
\cite[Section 5]{Tanaka-Tasaki},
Theorem \ref{thm:main1} yields

\begin{thm} \label{thm:main2}
Let $M$ be an irreducible Hermitian symmetric space of compact type
and $L_0, L_1$ be real forms of $M$ which intersect transversally.
Then the following results hold.
\begin{enumerate}
\item If $M=G_{2m}^{\mathbb C}(\mathbb C^{4m}) (m \geq 2)$,
$L_0$ is congruent to $G_{m}^{\mathbb H}(\mathbb H^{2m})$ and
$L_1$ is congruent to $U(2m)$, then we have
$$
HF(L_0, L_1 : {\mathbb Z}_2) \cong ({\mathbb Z}_2)^{2^m},
$$
where $2^m < {2m \choose m} =\#_2L_0 < 2^{2m} =\#_2L_1$.
Here $\#_2L$ denotes the 2-number of $L$.
\item Otherwise, we have
$$
HF(L_0, L_1 :{\mathbb Z}_2) \cong ({\mathbb Z}_2)^{\min\{ \#_2L_0, \#_2L_1 \}}.
$$
\end{enumerate}
\end{thm}

In the case where $M$ is non-irreducible,
although we must determine the rank of
$HF(L_0, L_1 :{\mathbb Z}_2)$ case by case (see Section 5),
at least we can state in general the following fact by combining
Theorem \ref{thm:main1} with Lemma 3.1 in \cite{Tasaki10}.

\begin{cor}
Let $M$ be a monotone Hermitian symmetric space of compact type and
$L_0, L_1$ be real forms of $M$ whose minimal Maslov numbers are greater
than or equal to $3$. Then for any $\phi \in \mathrm{Ham}(M,\omega)$ we have
\begin{equation*}
L_0 \cap \phi L_1 \neq \emptyset.
\end{equation*}
\end{cor}

Theorems \ref{thm:main1} and \ref{thm:main2} can be regarded as
a solution for a problem proposed by Y.-G. Oh in \cite[Section 6]{Oh93}.
Here we review the definition of 2-number introduced
by Chen and Nagano \cite{Chen-Nagano88}.
A subset $S$ in a Riemannian symmetric space $M$ is called an
{\it antipodal set}, if the geodesic symmetry $s_x$ fixes every point
of $S$ for any point $x$ of $S$.
The {\it 2-number} $\#_2M$ of $M$ is defined as the supremum of
the cardinalities of antipodal sets in $M$, which is known to be finite.
An antipodal set in $M$ is said to be {\it great} if its cardinality
attains $\#_2M$.
Takeuchi \cite{Takeuchi89} proved that if $L$ is a symmetric $R$-space,
then
\begin{equation*}
\#_2L = SB(L,\mathbb Z_2)
\end{equation*}
holds.
Note that any real form of Hermitian symmetric spaces of compact type
is a symmetric $R$-space, which is shown in \cite{Takeuchi84}. 
These facts and the invariance of $HF(L_0, L_1 : {\mathbb Z}_2)$
under Hamiltonian isotopies of $M$ imply

\begin{cor} \label{cor:GAG}
Let $M$ be an irreducible Hermitian symmetric space of compact type
and $(L_0, L_1)$ be a pair of real forms of $M$.
Then for any $\phi \in \mathrm{Ham}(M,\omega)$ such that
$L_0$ and $\phi L_1$ intersect transversally,
the following inequalities hold.
\begin{enumerate}
\item If $M=G_{2m}^{\mathbb C}(\mathbb C^{4m}) (m \geq 2)$,
$L_0$ is congruent to $G_{m}^{\mathbb H}(\mathbb H^{2m})$ and
$L_1$ is congruent to $U(2m)$, then we have
$$
\#(L_0 \cap \phi L_1) \geq 2^m.
$$
\item Otherwise, we obtain
\begin{equation} \label{eq:GAG}
\#(L_0 \cap \phi L_1) \geq \min\{ SB(L_0,\mathbb Z_2), SB(L_1,\mathbb Z_2) \}.
\end{equation}
\end{enumerate}
\end{cor}

As the case (1) above shows, in general, we cannot estimate
$\#(L_0 \cap \phi L_1)$
by the sum of ${\mathbb Z}_2$-Betti numbers of $L_0$ or $L_1$.
The estimate in the case (1) is sharp.
Note that we can construct many examples which
do not satisfy inequality (\ref{eq:GAG}) for the reducible case
(see Section 5).
We call inequality (\ref{eq:GAG})
{\it the generalized Arnold-Givental inequality}.
Indeed, (\ref{eq:GAG}) yields

\begin{cor}[\bf Oh \cite{Oh95} and \cite{Oh93} Theorem 1.3] \label{cor:Oh}
Let $(M,J_0,\omega)$ be an irreducible Hermitian symmetric space of
compact type and $\sigma$ be an anti-holomorphic involutive isometry of $M$.
Then Conjecture \ref{conj:AG} is true for the real form
$L=\mathrm{Fix}(\sigma)$ of $M$.
\end{cor}

\begin{rem} \rm
Real forms of Hermitian symmetric spaces of compact type $M$
are classified by Leung \cite{Leung79} and Takeuchi \cite{Takeuchi84}.
If $M$ is irreducible, then real forms $L_0$ and $L_1$ of $M$,
which are not congruent each other, are given in the list below.
Hence we can apply Theorem \ref{thm:main2} and Corollary \ref{cor:GAG}
to the following.
$$
\begin{array}{|c|c|c|c|}
\hline
M & L_0 & L_1 & \#(L_0 \cap L_1) \\
\hline
G_{2q}^{\mathbb C}(\mathbb C^{2m+2q})
& G_q^{\mathbb H}(\mathbb H^{m+q})
& G_{2q}^{\mathbb R}(\mathbb R^{2m+2q})
& {m+q \choose q} \\
\hline
G_n^{\mathbb C}(\mathbb C^{2n})
& U(n) & G_n^{\mathbb R}(\mathbb R^{2n})
& 2^n \\
\hline
G_{2m}^{\mathbb C}(\mathbb C^{4m})
& G_m^{\mathbb H}(\mathbb H^{2m}) & U(2m)
& 2^m \\
\hline
SO(4m)/U(2m) & U(2m)/Sp(m) & SO(2m)
& 2^m \\
\hline
Sp(2m)/U(2m) & Sp(m) & U(2m)/O(2m)
& 2^m \\
\hline
Q_n(\mathbb C) & S^{k,n-k} & S^{l,n-l}
& 2k+2 \ (\mathrm{if} \ k \leq l) \\
\hline
E_6/T\cdot Spin(10) & F_4/Spin(9) & G_2^{\mathbb H}(\mathbb H^4)/\mathbb Z_2
& 3 \\
\hline 
E_7/T\cdot E_6 & T\cdot (E_6/F_4) & (SU(8)/Sp(4))/\mathbb Z_2
& 8 \\
\hline
\end{array}
$$
Here, $G_r^{\mathbb K}(\mathbb K^{n+r})$ denotes the Grassmann
manifold of $r$-planes in ${\mathbb K}^{n+r}$ over the field
${\mathbb K}={\mathbb R}$, $\mathbb C$ or $\mathbb H$.
We denote the $n$-dimensional complex hyperquadric and a real form of it
by $Q_n(\mathbb C)$ and
$\displaystyle S^{k,n-k}=(S^{k} \times S^{n-k})/\mathbb Z_2$,
respectively (see \cite{Tasaki10}).
\end{rem}

This paper is organized as follows.
The calculation of $HF(L_0,L_1:{\mathbb Z}_2)$ is based on
the Floer homology theory for monotone Lagrangian submanifolds
as developed in \cite{Oh93}.
Section 2 contains an overview of the Lagrangian Floer theory.
In Section 3, we give a criterion for a Hermitian symmetric
space of compact type $M$ to be monotone
(see Proposition \ref{pro:monotone}).
Using it, we can also calculate some examples in the case where
$M$ is {\it non-irreducible}.
They are exhibited in Section 5.
In Section 4, we shall prove Theorem \ref{thm:main1}.
There, we see that a geodesic symmetry of a Hermitian symmetric space $M$
induces a free $\mathbb Z_2$-action
on the space of $J_0$-holomorphic strips.
In the last section,
as an application of inequality (\ref{eq:GAG}),
we prove that a real form $S^{0,n}$ of $Q_n(\mathbb C)$ is globally
volume minimizing under Hamiltonian deformations.

\section{Lagrangian Floer Homology}

In this section, we briefly review the Lagrangian Floer theory
as developed in \cite{Oh93} (see also \cite{Al08}).
Let $(M,\omega)$ be a closed symplectic manifold, $L_0$ and $L_1$
two closed Lagrangian submanifolds which intersect transversally.
We choose a time-dependent family $J= \{ J_t \}_{0 \leq t \leq 1}$  
of almost complex structures on $M$ compatible with
the symplectic form $\omega$.
The Floer chain complex $CF(L_0,L_1)$ is the vector space over
$\mathbb Z_2$ generated by the finitely many elements of $L_0 \cap L_1$.
A {\it $J$-holomorphic strip} is a map
$u:{\mathbb R} \times [0,1] \to M$
which satisfies the equation
\begin{equation} \label{eq:strip1}
 \bar{\partial}_J u :=
  \frac{\partial u}{\partial s} + J_t(u) \frac{\partial u}{\partial t} = 0
\end{equation}
with the following boundary conditions:
\begin{eqnarray}
 & & u(\cdot,0) \in L_0,\ u(\cdot,1) \in L_1,\label{eq:strip2} \\
 & & u(-\infty,\cdot), u(+\infty,\cdot) \in L_0 \cap L_1,\label{eq:strip3}
\end{eqnarray}
where ${\mathbb R} \times [0,1]$ is regarded as a subset of $\mathbb C$
with coordinates $s+\sqrt{-1}t$.
A solution of the equation (\ref{eq:strip1}) with Lagrangian
boundary condition (\ref{eq:strip2}) satisfies the asymptotic condition
(\ref{eq:strip3}) if and only if the energy of $u$ 
\begin{equation*}
E(u) = \frac{1}{2} \int_{{\mathbb R} \times [0,1]}
 \left( \left|\frac{\partial u}{\partial s}\right|^2
  +\left|\frac{\partial u}{\partial t}\right|^2 \right)
\end{equation*}
is finite.
The space of all $J$-holomorphic strips that connect $p \in L_0 \cap L_1$ to
$q \in L_0 \cap L_1$ is denoted by $\tilde{\mathcal M}_J(L_0,L_1:p,q)$.
We set
\begin{equation*}
\tilde{\mathcal M}_J(L_0,L_1)
 := \bigcup_{p,q \in L_0 \cap L_1} \tilde{\mathcal M}_J(L_0,L_1:p,q).
\end{equation*}
A family of almost complex structures $J$ is said to be {\it regular}
if the linearization $D_u \bar{\partial}_J$ of $\bar{\partial}_J$
is surjective for all $u \in \tilde{\mathcal M}_J(L_0,L_1)$.
For a regular $J$, each $\tilde{\mathcal M}_J(L_0,L_1:p,q)$
is a finite-dimensional smooth manifold,
with connected components of different dimensions.
We denote by ${\mathcal J}^{reg}$
the set of all regular almost complex structures on $M$.
The set ${\mathcal J}^{reg}$ is a set of the second category
in the set of families of almost complex structures on $M$.
From now on we assume $J \in {\mathcal J}^{reg}$.
If $u \in \tilde{\mathcal M}_J(L_0,L_1:p,q)$, then
\begin{equation*}
\dim(T_u \tilde{\mathcal M}_J(L_0,L_1:p,q))
 = \mathrm{Index}(D_u \bar{\partial}_J).
\end{equation*}
The right hand side is the index of $D_u \bar{\partial}_J$.
It is the spectral flow of $\bar{\partial}_J$ along $u$
and is equal to the Maslov index $\mu(u)$ of $u$.

The moduli space $\tilde{\mathcal M}_J(L_0,L_1)$ has a natural action
of $\mathbb R$ by translation in the first variable.
Hence, we define
\begin{eqnarray*}
{\mathcal M}_J(L_0,L_1:p,q)
 &:=& \tilde{\mathcal M}_J(L_0,L_1:p,q)/{\mathbb R},\\
{\mathcal M}_J(L_0,L_1)
 &:=& \tilde{\mathcal M}_J(L_0,L_1)/{\mathbb R}.
\end{eqnarray*}
An {\it isolated trajectory} is a trajectory $u$ in
$\tilde{\mathcal M}_J(L_0,L_1)$ such that the equivalence class $[u]$
is a 0-dimensional component of ${\mathcal M}_J(L_0,L_1)$.
The boundary operator $\partial : CF(L_0,L_1) \to CF(L_0,L_1)$
is defined by
\begin{equation*}
 \partial(p) = \sum_{q \in L_0 \cap L_1} n(p,q) \cdot q,
\end{equation*}
where $n(p,q)$ is the mod-$2$ number of isolated trajectories in
$\tilde{\mathcal M}_J(L_0,L_1:p,q)$.

To define the Floer homology group with coefficients in ${\mathbb Z}_2$
\begin{equation*}
HF(L_0,L_1:{\mathbb Z}_2)
 := \frac{\mathrm{Ker}(\partial)}{\mathrm{Im}(\partial)},
\end{equation*}
we must assume some topological conditions on $M, L_0$ and $L_1$.
For a closed Lagrangian submanifold $L$ in a symplectic manifold
$(M,\omega)$, two homomorphisms
\begin{equation*}
I_{\mu,L} : \pi_2(M,L) \to {\mathbb Z},\ \ 
I_\omega : \pi_2(M,L) \to {\mathbb R}
\end{equation*}
are defined as follows.
For a smooth map $w : (D^2,\partial D^2) \to (M,L)$,
$I_{\mu,L}(w)$ is defined to be the Maslov number of the bundle pair
$(w^*TM,(w|\partial D^2)^*TL)$
and $I_\omega$ is defined by $I_\omega(w)=\int_{D^2}w^*\omega$.
Then $L$ is said to be {\it monotone} if there exists a constant $\alpha>0$
such that $I_\omega = \alpha I_{\mu,L}$.
The {\it minimal Maslov number} $\Sigma_L$ of $L$ is defined to be
the positive generator of $\mathrm{Im}(I_{\mu,L}) \subset {\mathbb Z}$.
Oh proved the following

\begin{thm}[\cite{Oh93} Theorems 4.4, 5.1] \label{thm:Oh}
Let $(L_0,L_1)$ be a pair of monotone Lagrangian submanifolds
which intersect transversally.
Suppose that $\Sigma_{L_i} \geq 3$ for $i=0,1$ and
$\mathrm{Im}(\pi_1(L_i)) \subset \pi_1(M)$ is a torsion subgroup
for at least one of $i=0,1$.
Then there exists a dense subset ${\mathcal J}' \subset {\mathcal J}^{reg}$
such that if $J \in {\mathcal J}'$, then we have

\rm{(1)}\ $\partial$ is well-defined,

\rm{(2)}\ $\partial \circ \partial=0$,

\rm{(3)}\ $HF(L_0,L_1:{\mathbb Z}_2)$ is independent of $J$ and
Hamiltonian isotopies.
\end{thm}

Let $M$ be a Hermitian symmetric space of compact type.
Since $M$ is simply connected,
the condition that $\mathrm{Im}(\pi_1(L_i)) \subset \pi_1(M)$
is automatically satisfied.
Therefore, to apply Theorem \ref{thm:Oh} to real forms $L_0, L_1$ of $M$,
it suffices to assume that $L_0$ and $L_1$ are monotone and
$\Sigma_{L_i} \geq 3$ for $i=0,1$.
Moreover, we can specify the case where a real form $L$ does not satisfy
the condition that $\Sigma_L \geq 3$ from arguments in Section 3.
If $M$ is irreducible, the only exceptional case is $L={\mathbb R}P^1$
in $M={\mathbb C}P^1$, where $\Sigma_L=2$.
Hence, a real form $L$ of a Hermitian symmetric space $M$ of compact type
does not satisfy the condition that $\Sigma_L \geq 3$ if and only if
$M$ has ${\mathbb C}P^1$ as an irreducible factor and 
${\mathbb R}P^1 \subset {\mathbb C}P^1$ is an irreducible factor of
the real form $L$.

\begin{rem} \rm
If the assumption that $\Sigma_{L_i} \geq 3$ is not satisfied,
then we have to analyze the structure of disc bubbles to prove that
$\partial \circ \partial=0$.
It requires the classification of holomorphic discs with
Maslov index $2$ (see \cite{Al08}).
\end{rem}

\section{Monotonicity and minimal Maslov number of a real form}

In this section,
let $(M,J,\omega)$ be a compact K\"ahler manifold with
complex structure $J$ and K\"ahler form $\omega$.
The first Chern class of $(M,J,\omega)$ is denoted by
$c_1(M):=c_1(TM,J)$.
Then two homomorphisms
\begin{equation*}
I_c : \pi_2(M) \to {\mathbb Z},\ \ 
I_\omega : \pi_2(M) \to {\mathbb R}
\end{equation*}
are defined as follows.
For a smooth map $u : S^2 \to M$ which represents an element
$A \in \pi_2(M)$,
$I_c(A)$ is defined to be the Chern number
$c_1(A):=\langle c_1(M),[u] \rangle$
and $I_\omega$ is defined by $I_\omega(A)=\int_{S^2}u^*\omega$
as in the case of Lagrangian submanifolds.
Then $(M,J,\omega)$ is said to be {\it monotone} if
there exists a positive constant $\alpha>0$ such that
$I_\omega = \alpha I_c$.
The {\it minimal Chern number} $\Gamma_{c_1}$ of $M$ is defined to be
the positive generator of the subgroup $I_c(\pi_2(M))$ of $\mathbb Z$.

The Ricci form $\rho$ of $(M,J,\omega)$
is a closed $(1,1)$-form on $M$ and $\rho/2\pi$ represents
the first Chern class $c_1(M) \in H^2(M,{\mathbb Z})$.
$(M,J,\omega)$ is called {\it K\"ahler-Einstein} if there exists
a constant $c$ such that $\rho=c\omega$.
It is straightforward to check that a K\"ahler-Einstein manifold
$(M,J,\omega)$ with a positive Ricci constant $c$ is monotone.
In particular, an {\it irreducible} Hermitian symmetric space
of compact type is monotone.

Now we shall give a useful criterion for a Hermitian symmetric space
of compact type to be monotone.
A Hermitian symmetric space of compact type $(M,J_0,\omega)$
can be decomposed as 
$$
(M,J_0,\omega)
\cong (M_1,J_1,\omega_1) \times (M_2,J_2,\omega_2)
\times \cdots \times (M_k,J_k,\omega_k),
$$
where each $(M_i,J_i,\omega_i)$ is an irreducible one.
Then the K\"ahler form $\omega$ and the Ricci form $\rho$ of $M$ are
represented as
$\omega = \omega_1 \oplus \omega_2 \oplus \cdots \oplus \omega_k$ and
$\rho = \rho_1 \oplus \rho_2 \oplus \cdots \oplus \rho_k$,
respectively.

\begin{pro} \label{pro:monotone}
Let $(M,J_0,\omega)$ be a Hermitian symmetric space of compact type.
Then $(M,\omega)$ is monotone as a symplectic manifold if and only if
$(M,J_0,\omega)$ is a K\"ahler-Einstein manifold with
a positive Ricci constant.
\end{pro}

\noindent
{\bf Proof.}\quad
Since each irreducible component $M_i$ of $M$ is K\"ahler-Einstein,
there exist constants $c_i>0$ such that
$\rho_i = c_i \omega_i$ for $i=1,2,\ldots,k$.
Then we have
\begin{eqnarray*}
[\omega]
= [\omega_1]+[\omega_2]+\cdots+[\omega_k]
= \frac{2\pi}{c_1}c_1(M_1)+\frac{2\pi}{c_2}c_1(M_2)
 +\cdots+\frac{2\pi}{c_k}c_1(M_k).
\end{eqnarray*}
If $M$ is K\"ahler-Einstein, i.e., $c_1=\cdots=c_k=:c$, then
$[\omega]=(2\pi/c) c_1(M)$ holds.
It shows that $M$ is monotone.

Conversely, if $M$ is monotone, then there exists a constant $\alpha>0$
such that
$$
\frac{2\pi}{c_1}c_1(M_1)+\frac{2\pi}{c_2}c_1(M_2)
 +\cdots+\frac{2\pi}{c_k}c_1(M_k)
= \alpha(c_1(M_1)+c_1(M_2)+\cdots+c_1(M_k)).
$$
It yields that $2\pi/c_i=\alpha$ for $i=1,2,\ldots,k$,
and hence $M$ is a K\"ahler-Einstein manifold with positive Ricci constant
$2\pi/\alpha$.
\hfill \qed

\smallskip

The following formula is necessary for us to ensure the monotonicity for
a real form $L$ and to estimate its minimal Maslov number $\Sigma_L$
(see \cite[Lemma 2.1]{Oh93}).

\begin{lem}[Viterbo] \label{lem:Viterbo}
Let $(M,J,\omega)$ be a compact K\"ahler manifold
and $L$ a closed Lagrangian submanifold.
Let $w,w' : (D^2,\partial D^2) \to (M,L)$ be smooth maps of pairs
satisfying
$w|_{\partial D^2}=w'|_{\partial D^2}$.
If we define a map $u$ from $S^2=D^2 \cup \overline{D^2}$ to $M$ as
$$
u(z)=
\left\{
 \begin{array}{rrr}
  w(z), & & z \in D^2,\\
  w'(z), & & z \in \overline{D^2},
 \end{array}
\right.
$$
then we have
$$
I_{\mu,L}(w)-I_{\mu,L}(w')=2c_1([u]).
$$
\end{lem}

\begin{cor}
Let $(M,J,\omega)$ be a monotone compact K\"ahler manifold.
Then the fixed point set $L=\mathrm{Fix}(\sigma)$
of an involutive anti-holomorphic isometry $\sigma$ is monotone.
\end{cor}

\noindent
{\bf Proof.}\quad
For any $A \in \pi_2(M,L)$, we take a smooth map
$w : (D^2,\partial D^2) \to (M,L)$
as a representative of $A$.
Then we can define another smooth map
$w'=\sigma \circ w : (D^2,\partial D^2) \to (M,L)$.
By Lemma \ref{lem:Viterbo}, we have
\begin{eqnarray} \label{eq:MaslovChern}
I_{\mu,L}(w)=c_1([u]).
\end{eqnarray}
Since $M$ is monotone, there exists a constant $\alpha>0$ such that
$\int_{S^2}u^*\omega = \alpha c_1([u])$, and hence
$\int_{S^2}u^*\omega=\alpha I_{\mu,L}(w)$.
The left hand side of this equation is equal to $2I_\omega(A)$.
Therefore, $I_\omega(A)=\frac{\alpha}{2} I_{\mu,L}(w)$.
That is, $L$ is a monotone Lagrangian submanifold with the monotone
constant $\frac{\alpha}{2}$.
\hfill \qed

\smallskip

The definitions of minimal Maslov and Chern numbers and
equality (\ref{eq:MaslovChern}) immediately imply

\begin{cor}
For a compact K\"ahler manifold $(M,J,\omega)$,
the minimal Chern number $\Gamma_{c_1}$ of $M$ and the minimal Maslov
number $\Sigma_L$ of a real form $L$ of $M$ satisfy
$$
\Sigma_L \geq \Gamma_{c_1}.
$$
\end{cor}

Therefore, to apply Theorem \ref{thm:Oh}
to a pair of real forms $(L_0,L_1)$ of
a Hermitian symmetric space $M$ of compact type,
it suffices to {\it assume that $M$ is monotone}
(it is equivalent that $M$ is K\"ahler-Einstein)
{\it and each $L_i$}
(it is automatically monotone Lagrangian submanifold)
{\it satisfies $\Sigma_{L_i} \geq 3$}.

The minimal Chern numbers of irreducible Hermitian symmetric spaces
$M$ of compact type are calculated as follows
(see \cite[p. 521]{Borel-Hirzebruch58}).

$$
\begin{array}{|c|c|c|}
\hline
M & \Gamma_{c_1} \\
\hline
U(m+n)/(U(m) \times U(n)) & m+n \\
\hline
SO(2n)/U(n) & 2n-2 \\
\hline
Sp(n)/U(n) & n+1\\
\hline
SO(n+2)/(SO(2) \times SO(n)) & n \\
\hline
E_6/T\cdot Spin(10) & 12 \\
\hline
E_7/T\cdot E_6 & 18 \\
\hline
\end{array}
$$
Therefore, any real form $L$ of $M$ satisfies that
$\Sigma_L \geq 3$ except for $L={\mathbb R}P^1$ in
$M={\mathbb C}P^1=U(2)/(U(1) \times U(1))$.
This case is treated in \cite[Section 5]{Oh93'} independently.

\section{Calculation of the Floer homology}

In this section,
we consider a monotone Hermitian symmetric space $(M,J_0,\omega)$
of compact type with the standard complex structure $J_0$ and
the standard K\"ahler form $\omega$.
Let $L_0$ and $L_1$ be real forms of $M$ which intersect transversally
and satisfy that $\Sigma_{L_i} \geq 3$ for $i=0,1$.
We take $J_0=J_t$ for all $t \in [0,1]$.
The following result ensures that $J_0$ can be used to calculate
the Floer homology $HF(L_0, L_1 : {\mathbb Z}_2)$
\cite[Main Theorem]{Oh97}.

\begin{thm}[Regularity \cite{Oh97}] \label{thm:Oh97}
Let $(M,J_0,\omega)$ be a K\"ahler manifold with non-negative
holomorphic bisectional curvature.
Let $L_0$ and $L_1$ be closed totally geodesic Lagrangian submanifolds
in $M$ which intersect transversally.
Then the complex structure $J_0$ is regular, i.e.,
the linearization $D_u \bar{\partial}_{J_0}$ of $\bar{\partial}_{J_0}$
is surjective for all $u \in \tilde{\mathcal M}_{J_0}(L_0,L_1)$.
\end{thm}

We apply the above theorem to the case where $(M,J_0,\omega)$
is a Hermitian symmetric space of compact type.
By the same argument as \cite[Proposition 4.5]{Oh95}
for the case where $(L_0,L_1)=(L,\phi(L))$,
Theorem \ref{thm:Oh97} yields

\begin{pro}[Compactness] \label{pro:compactness}
Let $(M,J_0,\omega)$ be a monotone Hermitian symmetric space
of compact type.
Let $L_0$ and $L_1$ be real forms of $M$ which intersect transversally.
In addition, assume that $\Sigma_{L_i} \geq 3$ for $i=0,1$.
Then the $0$-dimensional part of ${\mathcal M}_{J_0}(L_0,L_1)$ is compact
and the $1$-dimensional part of ${\mathcal M}_{J_0}(L_0,L_1)$
is compact up to the splitting of two isolated trajectories.
Therefore, $\partial_{J_0}^2=0$.
\end{pro}

The following Theorem by Tanaka and Tasaki is essential for calculation.

\begin{thm}[Theorem 1.1 in \cite{Tanaka-Tasaki}] \label{thm:Tanaka-Tasaki}
Let $M$ be a Hermitian symmetric space of compact type and
$L_0$ and $L_1$ real forms which intersect transversally.
Then the intersection $L_0 \cap L_1$ is an antipodal set of $L_0$ and $L_1$.
\end{thm}

The geodesic symmetry $s_p$ at any point $p$
of a Hermitian symmetric space is a holomorphic isometry.
In Theorem \ref{thm:Tanaka-Tasaki} the intersection $L_0 \cap L_1$
is also an antipodal set in $M$, because $L_0$ and $L_1$ are
totally geodesic, which yields the following:

\begin{lem} \label{lem:z_2action}
Under the assumption of Theorem \ref{thm:Tanaka-Tasaki},
for any $p \in L_0 \cap L_1$,
where $L_0 \cap L_1$ is not empty by Lemma 3.1 in \cite{Tasaki10},
the geodesic symmetry $s_p$ satisfies
$$
s_p(L_0) = L_0, \quad s_p(L_1) = L_1, \quad
s_p(q) = q \quad (q \in L_0 \cap L_1).
$$
\end{lem}

\noindent
{\bf Proof.}\quad
Since a real form of $M$ is totally geodesic,
we have $L_i = \mathrm{Exp}_p(T_pL_i)$ for $i=0,1$.
Remark that $s_p^2=\mathrm{id_M}$.
Since the differential map of $s_p$ satisfies $(ds_p)_p = -1$,
we obtain
$$
s_p(L_i) = \mathrm{Exp}_p((ds_p)_pT_pL_i)
= \mathrm{Exp}_p(T_pL_i) = L_i.
$$
By Theorem \ref{thm:Tanaka-Tasaki},
the intersection $L_0 \cap L_1$ is an antipodal set of $L_0$ and $L_1$.
Hence, by definition, $s_xy = y$ holds for any $x, y \in L_0 \cap L_1$.
In particular, we have $s_p(q) = q$ for $q \in L_0 \cap L_1$.
\hfill \qed

\smallskip

Now we shall calculate $HF(L_0, L_1 : {\mathbb Z}_2)$.
By assumption, the intersection $L_0 \cap L_1$ consists of finite points.
We choose any two points $p,q \in L_0 \cap L_1$.
By Lemma \ref{lem:z_2action},
we see that $p,q$ are fixed points of the action of $s_p$.
Let $u$ be a $J_0$-holomorphic strip in 
$\tilde{\mathcal M}_{J_0}(L_0,L_1:p,q)$.
It satisfies the boundary conditions
\begin{equation*}
u(s,0) \in L_0,\ u(s,1) \in L_1,\ u(-\infty,t)=p,\ u(+\infty,t)=q.
\end{equation*}
Using the holomorphic isometry $s_p$,
let us define another holomorphic map
$\bar{u}(s,t):=s_p(u(s,t))$ from ${\mathbb R} \times [0,1]$ to $M$.
By Lemma \ref{lem:z_2action},
real forms $L_0$ and $L_1$ are invariant under the action of $s_p$.
Hence, the holomorphic map $\bar{u}$ satisfies that
\begin{equation*}
\bar{u}(s,0) = s_p(u(s,0)) \in L_0,\
\bar{u}(s,1) = s_p(u(s,1)) \in L_1
\end{equation*}
and $\bar{u}(-\infty,t) = s_p(u(-\infty,t))=s_p(p)=p,\ \bar{u}(+\infty,t) = s_p(u(+\infty,t))=s_p(q)=q$.
It says that the holomorphic map $\bar{u}$ also belongs to
$\tilde{\mathcal M}_{J_0}(L_0,L_1:p,q)$.
Moreover, $s_p \circ \bar{u}=u$ and we see that
$[\bar{u}] \neq [u] \in {\mathcal M}_{J_0}(L_0,L_1:p,q)$
from the definition of the map $s_p$,
and hence the moduli space ${\mathcal M}_{J_0}(L_0,L_1:p,q)$
possesses a free $\mathbb Z_2$-action induced from $s_p$.
Since the $0$-dimensional part of the moduli space
${\mathcal M}_{J_0}(L_0,L_1:p,q)$ is compact
by Proposition \ref{pro:compactness},
it contains an even number of elements.
Therefore, we obtain

\begin{pro}[Vanishing] \label{pro:vanishing}
Under the same assumptions as in Proposition \ref{pro:compactness},
the number of $0$-dimensional components of
${\mathcal M}_{J_0}(L_0,L_1:p,q)$ are even
and so the boundary operator
$\partial : CF(L_0,L_1) \to CF(L_0,L_1)$
vanishes.
\end{pro}

Thus we complete the proof of Theorem \ref{thm:main1}.

\begin{rem} \rm
If real forms $L_0$ and $L_1$ are congruent, then the above calculation
provides us with an alternative proof of the known
fact that
$$
HF(L_0, L_0 :{\mathbb Z}_2) \cong ({\mathbb Z}_2)^{\#_2L_0}
 = ({\mathbb Z}_2)^{SB(L_0, \mathbb Z_2)},
$$
because the intersection $L_0 \cap L_1$ is a great antipodal set of
$L_0$ (and $L_1$), which is proved in \cite[Theorem 1.3]{Tanaka-Tasaki}
\end{rem}

\section{Some examples for the reducible case}

Let $(M,J_0,\omega)$ be an irreducible Hermitian symmetric space
of compact type and $\sigma : M \to M$ an involutive anti-holomorphic 
isometry.
Since the product $M\times M$ of $M$ is a K\"ahler-Einstein manifold
with positive Ricci constant,
we can apply Theorem \ref{thm:main1} to real forms of $M\times M$.
Since $(x, y) \mapsto (\sigma(y), \sigma(x))$
is an involutive anti-holomorphic isometry of $M \times M$,
whose fixed point set
$$
D_\sigma(M) = \{(x, \sigma(x)) \mid x \in M\}
$$
is a real form of $M\times M$.
On the other hand, for real forms $L_0$ and $L_1$ of $M$,
we see that $L_0 \times L_1$ is a real form of $M \times M$.
Then
$$
(L_0 \times L_1) \cap D_\sigma(M)
= \{(x, \sigma(x)) \mid x \in L_0 \cap \sigma^{-1}(L_1)\}.
$$
The condition that two real forms $L_0 \times L_1$ and $D_\sigma(M)$
of $M \times M$ intersect transversally
is equivalent to
the condition that two real forms $L_0$ and $\sigma^{-1}(L_1)$
of $M$ intersect transversally.
In this situation, we obtain
$$
\#\{(L_0 \times L_1) \cap D_\sigma(M)\}
= \#\{L_0 \cap \sigma^{-1}(L_1)\}.
$$
Moreover, $\sigma^{-1}(L_1)$ is congruent to $L_1$.

\begin{exam} \rm
Let $M$ be the complex projective space ${\mathbb C}P^n$.
Real forms $L_0$ and $L_1$ of $M$ are congruent to ${\mathbb R}P^n$.
Then
$
\#\{(L_0 \times L_1) \cap D_\sigma(M)\}
= \#\{L_0 \cap \sigma^{-1}(L_1)\}
= n + 1
$.
By Lemma 1.1 in \cite{Chen-Nagano88}, we have
$$
\#_2(L_0 \times L_1) = \#_2(L_0) \#_2(L_1) = (n + 1)^2, \quad
\#_2(D_\sigma(M)) = \#_2 M = n + 1,
$$
and hence the intersection number of the two real forms is equal to
smaller 2-number $n+1$.
Moreover, we can easily check that $I_{\mu,D_\sigma(M)}=2(n+1)$
and $I_{\mu,L_0 \times L_1} \geq 3$ for $n \geq 2$.
By Theorem \ref{thm:main1}, we have
$$
HF(L_0 \times L_1, D_\sigma(M) : {\mathbb Z}_2) \cong ({\mathbb Z}_2)^{n+1}
$$
for $n \geq 2$.
When $n=1$, two real forms
$L_0 \times L_1 = {\mathbb R}P^1 \times {\mathbb R}P^1 \cong T^2$
and
$D_\sigma(M) \cong S^2$
can be regarded as real forms of $2$-dimensional complex hyperquadric
$Q_2(\mathbb C) \cong {\mathbb C}P^1 \times {\mathbb C}P^1$.
Although $\Sigma_{T^2}=2$, we can also prove that
$HF(S^2, T^2 : {\mathbb Z}_2) \cong {\mathbb Z}_2 \oplus {\mathbb Z}_2$
by combining the arguments in Section 4 and in \cite{Al08}.
Hence, the pair $(L_0 \times L_1,D_\sigma(M))$ satisfies the generalized
Arnold-Givental inequality (\ref{eq:GAG}).
\end{exam}

\begin{exam} \rm
Put $M = Q_n(\mathbb C)$.
Assume that real forms $L_0, L_1$ of $M$ are congruent to
$S^{k,n-k}, S^{l,n-l}\; (0 \le k \le l \le [n/2])$,
respectively.
Then by a result in \cite{Tasaki10}, we have
$$
\#\{(L_0 \times L_1) \cap D_\sigma(M)\}
= \#\{L_0 \cap \sigma^{-1}(L_1)\}
= 2(k + 1).
$$
If $n \geq 3$, then minimal Maslov numbers
$I_{\mu,D_\sigma(M)}$ and $I_{\mu,L_0 \times L_1}$
are greater than or equal to $3$.
By Theorem \ref{thm:main1}, we obtain
$$
HF(L_0 \times L_1, D_\sigma(M) : {\mathbb Z}_2) \cong ({\mathbb Z}_2)^{2(k+1)}.
$$
Note that $\#_2(L_0 \times L_1) = 4(k + 1)(l + 1)$ and
$\#_2(D_\tau(M)) = 2([n/2] + 1)$.
Hence the intersection number of the two real forms
coincides with $\min\{ \#_2(L_0 \times L_1), \#_2(D_\sigma(M)) \}$
only in the case where $k = l = [n/2]$,
otherwise the intersection number is smaller than it.

In this way, we can construct many pairs of real forms which
do not satisfy the generalized Arnold-Givental inequality (\ref{eq:GAG}).
\end{exam}

\section{A volume estimate for a real form under Hamiltonian isotopies}

In general, a closed Lagrangian submanifold $L$ in a K\"ahler manifold
$(M,J,\omega)$ is said to be {\it Hamiltonian volume minimizing}
if it satisfies
$$
\mathrm{vol}(\phi L) \geq \mathrm{vol}(L)
$$
for any Hamiltonian diffeomorphism $\phi \in \mathrm{Ham}(M,\omega)$
(see \cite{Oh90}).
Non-trivial known examples of Hamiltonian volume minimizing Lagrangian
submanifolds are very few.
It is known that the real form ${\mathbb R}P^n$ in the complex projective
space ${\mathbb C}P^n$ and real form $S^1 \times S^1$ in $S^2 \times S^2$
are Hamiltonian volume minimizing Lagrangian submanifolds
(see \cite{Oh90} and \cite{IOS2003}).
Since $S^2 \times S^2$ is isomorphic to $Q_2(\mathbb C)$,
it is worthwhile to verify which real form of the complex hyperquadric
$Q_n(\mathbb C)$ is Hamiltonian volume minimizing.
In fact, the Hamiltonian stabilities of real forms of $Q_n(\mathbb C)$
were determined by Oh \cite{Oh90}.

Here, we give a lower bound of the volume of the image $\phi(S^{k,n-k})$ of
a real form $\displaystyle S^{k,n-k}=(S^{k} \times S^{n-k})/\mathbb Z_2$ of
$Q_n(\mathbb C)$ by any $\phi \in \mathrm{Ham}(Q_n(\mathbb C),\omega)$.
By the generalized Arnold-Givental inequality (\ref{eq:GAG}), we have
\begin{equation} \label{eq:AGAQn}
\#(S^{0,n} \cap \phi S^{k,n-k}) \geq
 \min\{ SB(S^{0,n},\mathbb Z_2), SB(S^{k,n-k},\mathbb Z_2) \} =2.
\end{equation}
Moreover, we use the following Crofton type formula.

\begin{thm}[Le \cite{Le93}] \label{thm:Le93} \rm
Let $N$ be an $n$-dimensional real submanifold in 
$Q_n(\mathbb C) \cong \widetilde{G_n}(\mathbb R^{n+2})$.
Then
\begin{equation} \label{eq:Le}
\int_{SO(n+2)} \#(g S^n \cap N)d\mu_{SO(n+2)}(g)
 \leq 2 \frac{\mathrm{vol}(SO(n+2))}{\mathrm{vol}(S^n)} \mathrm{vol}(N)
\end{equation}
holds.
\end{thm}

\begin{pro} \label{pro:vol} \rm
For any $\phi \in \mathrm{Ham}(Q_n(\mathbb C),\omega)$, we have
$\mathrm{vol}(\phi S^{k,n-k}) \geq \mathrm{vol}(S^n)$.
\end{pro}

\noindent
{\bf Proof.}\quad
Put $N=\phi S^{k,n-k} \ (k=0,1,\ldots,[n/2])$.
Then (\ref{eq:Le}) and (\ref{eq:AGAQn}) yield
\begin{eqnarray*}
\mathrm{vol}(\phi S^{k,n-k})
 & \geq & \frac{\mathrm{vol}(S^n)}{2 \mathrm{vol}(SO(n+2))}
  \int_{SO(n+2)} \#(g S^n \cap \phi S^{k,n-k})d\mu_{SO(n+2)}(g) \\
 & \geq & \frac{\mathrm{vol}(S^n)}{2 \mathrm{vol}(SO(n+2))}
  \int_{SO(n+2)} 2 d\mu_{SO(n+2)}(g) \\
 &=& \mathrm{vol}(S^n).
\end{eqnarray*}
\hfill \qed

\smallskip

Gluck, Morgan and Ziller \cite{Gluck-Morgan-Ziller1989} proved that
$S^{0,n}=S^n$ in $Q_n(\mathbb C)$ is volume minimizing in its homology class
when $n$ is even.
On the other hand, since the homology $H_k(Q_n(\mathbb C))$ vanishes when $k$ is odd,
$S^{0,n}$ can not be homologically volume minimizing in $Q_n(\mathbb C)$
in the case where $n$ is odd.
At least, we can conclude from Proposition \ref{pro:vol} the following

\begin{cor}
A real form $S^{0,n}$ of the complex hyperquadric $Q_n(\mathbb C)$
is Hamiltonian volume minimizing.
\end{cor}

\section*{Acknowledgements}

We would like to thank Professor Yoshihiro Ohnita for pointing out
that Proposition \ref{pro:vol} produces new examples of Hamiltonian
volume minimizing Lagrangian submanifolds as in Corollary 24.
We would like to thank Professor Martin Guest for some helpful
comments on this paper.
The authors are also grateful to the referee, whose useful comments improved
the manuscript.
The first author was partly supported by
the Grant-in-Aid for Young Scientists (B) (No.~22740043),
MEXT.
The second author was partly supported by
the Grant-in-Aid for Young Scientists (B) (No.~23740057),
JSPS.
The third author was partly supported by
the Grant-in-Aid for Science Research (C) (No.~21540063),
JSPS.

\begin{flushleft}
H.~Iriyeh

{\sc School of Science and Technology for Future Life\\
Tokyo Denki University\\
Kanda-Nishiki-Cho, Chiyoda-Ku\\
Tokyo, 101-8457 Japan}

{\it e-mail} : {\tt hirie@im.dendai.ac.jp}
\end{flushleft}

\begin{flushleft}
T.~Sakai

{\sc Department of Mathematics and Information Sciences\\
Tokyo Metropolitan University\\
Minami-osawa, Hachioji-shi\\
Tokyo, 192-0397 Japan
}

{\it e-mail} : {\tt sakai-t@tmu.ac.jp}
\end{flushleft}

\begin{flushleft}
H.~Tasaki

{\sc Division of Mathematics\\
Faculty of Pure and Applied Sciences\\
University of Tsukuba\\
Tsukuba, Ibaraki, 305-8571 Japan}

{\it e-mail} : {\tt tasaki@math.tsukuba.ac.jp}
\end{flushleft}


\begin{thebibliography}{9}

\bibitem{Al08}
G. Alston,
{\em Lagrangian Floer homology of the Clifford torus and
real projective space in odd dimensions},
arXiv:0902.0197v2.

\bibitem{Alston-Amorim10}
G. Alston and L. Amorim,
{\em Floer cohomology of torus fibers and real Lagrangians
in Fano toric manifolds},
arXiv:1003.3651v1.

\bibitem{Arnold65}
V.I. Arnold,
{\em Sur une propri\'et\'e topologique des applications globalement
canoniques de la m\'ecanique classique},
C. R. Acad. Sc. Paris {\bf 261} (1965), 3719--3722.

\bibitem{Borel-Hirzebruch58}
A. Borel and F. Hirzebruch,
{\em Characteristic classes and homogeneous spaces I},
Amer. J. Math. {\bf 80} (1958), 458--538.

\bibitem{Chang-Jiang90}
K.-C. Chang and M.Y. Jiang,
{\em The Lagrange intersections for $({\mathbb C}P^n,{\mathbb R}P^n)$},
Manuscripta Math. {\bf 68} (1990), 89--100.

\bibitem{Chen-Nagano88}
B.-Y. Chen and T. Nagano,
{\em A Riemannian geometric invariant and its applications
to a problem of Borel and Serre},
Trans. Amer. Math. Soc. {\bf 308} (1988), 273--297.

\bibitem{Floer88}
A. Floer,
{\em Morse theory for Lagrangian intersections},
J. Differ. Geom. {\bf 28} (1988), 513--547.

\bibitem{FOOO}
K. Fukaya, Y.-G.Oh, H. Ohta and K. Ono,
{\em Floer theory of Lagrangian submanifolds over $\mathbb Z$},
preprint.

\bibitem{FOOO10}
K. Fukaya, Y.-G.Oh, H. Ohta and K. Ono,
{\em Toric degeneration and non-displaceable Lagrangian tori in
$S^2 \times S^2$},
arXiv:1002.1660v1.

\bibitem{Frauenfelder04}
U. Frauenfelder,
{\em The Arnold-Givental conjecture and moment Floer homology},
Int. Math. Res. Not. {\bf 42} (2004), 2179--2269.

\bibitem{Givental88}
A. Givental,
{\em Periodic maps in symplectic topology},
Funct. Anal. Appl. {\bf 23} (1989), 37--52.

\bibitem{Gluck-Morgan-Ziller1989}
H. Gluck, F. Morgan and W. Ziller,
{\em Calibrated geometries in Grassmann manifolds},
Comm. Math. Helv. {\bf 64} (1989), 256--268.

\bibitem{IOS2003}
H.~Iriyeh, H.~Ono and T.~Sakai,
{\em Integral geometry and Hamiltonian volume minimizing property of
a totally geodesic Lagrangian torus in $S^{2} \times S^{2}$},
Proc. Japan Acad. {\bf 79} Ser. A (2003), 167--170.

\bibitem{Leung79}
D.~P.~S.~Leung,
{\em Reflective submanifolds. IV, Classification of real forms
of Hermitian symmetric spaces},
J. Differ. Geom. {\bf 14} (1979), 179--185.

\bibitem{Le93}
L\^e H\^ong V\^an,
{\em Application of integral geometry to minimal surfaces},
Int. J. Math. {\bf 4} (1993), 89--111.

\bibitem{Oh90}
Y.-G. Oh,
{\em Second variation and stabilities of minimal
lagrangian submanifolds in K\"{a}hler manifolds},
Invent. Math. {\bf 101} (1990), 501-519.

\bibitem{Oh93}
Y.-G. Oh,
{\em Floer cohomology of Lagrangian intersections and pseudo-holomorphic disks, I},
Comm. Pure Appl. Math. {\bf 46} (1993), 949--993;
Addendum, Comm. Pure Appl. Math. {\bf 48} (1995), 1299--1302.

\bibitem{Oh93'}
Y.-G. Oh,
{\em Floer cohomology of Lagrangian intersections and pseudo-holomorphic disks, II:
$({\mathbb C}P^n, {\mathbb R}P^n)$},
Comm. Pure Appl. Math. {\bf 46} (1993), 995--1012.

\bibitem{Oh95}
Y.-G. Oh,
{\em Floer cohomology of Lagrangian intersections and pseudo-holomorphic disks, III:
Arnold-Givental conjecture},
The Floer Memorial volume, Progr. Math., vol. 133, Birkh\"auser, Basel 
(1995), 555--573.

\bibitem{Oh97}
Y.-G. Oh,
{\em Fredholm-Regularity of Floer's Holomorphic Trajectories on
K\"ahler Manifolds},
Kyungpook Math. J. {\bf 37} (1997), 153--164.

\bibitem{Takeuchi84}
M.~Takeuchi,
{\em Stability of certain minimal submanifolds of compact Hermitian
symmetric spaces},
Tohoku Math. J. {\bf 36} (1984), 293--314.

\bibitem{Takeuchi89}
M.~Takeuchi,
{\em Two-number of symmetric R-spaces},
Nagoya Math. J. {\bf 115} (1989), 43--46.

\bibitem{Tanaka-Tasaki}
M. S. Tanaka and H. Tasaki,
{\em The intersection of two real forms in Hermitian symmetric
spaces of compact type},
to appear in J. Math. Soc. of Japan.

\bibitem{Tasaki10}
H.~Tasaki,
{\em The intersection of two real forms in the complex
hyperquadric},
Tohoku Math. J. {\bf 62} (2010), 375--382.

\end{thebibliography}
\end{document}